\documentclass{amsart}
\usepackage{amsmath,amssymb,latexsym}
\usepackage[T1]{fontenc}
\usepackage{dsfont}
\usepackage{enumerate}

\usepackage{enumerate,xspace}

\usepackage{soul}

\usepackage{mdwlist}

\usepackage[pdftex]{hyperref}

\theoremstyle{plain}

\newtheorem{theoremA}{Theorem}

\newenvironment{claimproof}[1][Proof of Claim.] 
{%
	\proof[#1]%
	
}
{%
	\endproof%
}

\title{Finite groups contain large centralizers}
\author{Daniel Palac\'in}
\date{\today}
\address{ \, Abteilung f\"ur Mathematische Logik, Mathematisches Institut,
  Albert-Ludwig-Universit\"at Freiburg, Ernst-Zermelo-Stra\ss e 1, D-79104
  Freiburg, Germany}
\email{palacin@math.uni-freiburg.de}
\thanks{Research supported by MTM2017-86777-P as well as by the Deutsche
	Forschungsgemeinschaft (DFG, German Research Foundation) - 
	Project number 2100310301, part of the ANR-DFG 
	program GeoMod}
\begin{document}

\begin{abstract}
	Every finite non-abelian group of order $n$ has a non-central element whose centralizer has order exceeding $n^{1/3}$. The proof does not rely on the  classification of finite simple groups, yet it uses the Feit-Thompson theorem.
\end{abstract}

\maketitle

\section{Introduction}

A classical theorem of Brauer and Fowler \cite{BF55} states that a finite non-abelian group $G$ of even order with a center of odd order has a non-central element $x$ such that 
\[
|G|<|C_G(x)|^3.
\]
For finite non-abelian solvable groups, Bertram \cite{eB84} proved the same inequality  and asked whether the exponent $3$ could be improved to $2$. This question was answered affirmatively by Isaacs \cite{iI86}, who showed that every finite non-abelian solvable group contains a non-central element whose centralizer has order exceeding its index. 

In \cite{GR20}, Guralnick and Robinson considered some variants of the Brauer-Fowler theorem. Among other results, they prove in \cite[Theorem 5]{GR20} that any finite non-abelian group $G$ has a non-central element $x$ of $G$ such that
\[
|G| < \frac{6}{5}|C_G(x)|^3.
\]
Their proof does not rely on the classification of finite simple groups  but uses the Feit-Thompson odd order theorem as well as a degenerate case of a result of Griess \cite{rG78}. In fact, using the classification they slightly improved this result showing the following:
\begin{theoremA}\cite{GR20}\label{T:ThmGeneral}
	Let $G$ be a finite non-abelian group. Then, there exists a non-central element $x$ of $G$ such that
	\[
	|G|<|C_G(x)|^3.
	\]
\end{theoremA}
The purpose of this short note is to give a proof of this result without using the classification, but still using the Feit-Thompson theorem.

Note that as a consequence of the aforementioned result of Bertram \cite{eB84} (see also \cite[Lemma 5.1]{GR20}), to prove Theorem \ref{T:ThmGeneral} it suffices to consider finite non-solvable groups. Hence, since all finite groups of odd order are solvable by the Feit-Thompson odd order theorem, we are reduced to considering a finite non-abelian group $G$ such that $G/Z(G)$ has even order. Therefore, Theorem \ref{T:ThmGeneral} follows from the following statement:

\begin{theoremA}\label{T:ThmEven} 
Let $G$ be a finite non-abelian group of even order and let $t$ be a non-central element of $G$ such that $t^2$ is central. Then, there exists a non-central element $x$ of $G$ such that
\[
|G|\le |C_G(t)|^2\left( |C_G(x)| - \frac{1}{2} \right).
\]
\end{theoremA}

We remark that in general the exponents in Theorem \ref{T:ThmGeneral} and Theorem \ref{T:ThmEven} cannot be improved as occurs  in $\mathrm{SL}(2,2^n)$, where the centralizer of an involution has order $2^n$ and the maximum order of a centralizer of a non-identity element is $2^n+1$. 

\section{Proof of Theorem \ref{T:ThmEven}}

Let $G$ be a finite non-abelian group and let $t$ be a non-central element of $G$ such that $t^2$ is central. Write $Z=Z(G)$ and $k(G)$ for the number of conjugacy classes of $G$. Also, we let $i(Z)$ denote the number of involutions of $Z$. 

\subsubsection*{Claim}\label{Claim} The following equation holds:
\[
|G| \le (1 + i(Z))|C_G(t)| + (k(G) - |Z|) |C_G(t)|^2.
\] 
\begin{claimproof}
Let $W$ be the set of pairs $(x,y)$ in $G\times G$ such that $x$ is a conjugate of $t$ which inverts $y$, that is
\[
W =\left\{ (x,y)\in t^G \times G \, : \,  y^x = y^{-1} \right\},
\]
and set $W_y = \left\{ x \in t^G : (x,y)\in W\right\}$ for an element $y$ of $G$.

It is clear that a central element $y$ of $G$ with $W_y\neq \emptyset$ must be an involution or the identity element and in any case $W_y = t^G$. For an arbitrary involution $y$ of $G$, the set $W_y$ equals $t^G\cap C_G(y)$ and for any other element $y$ of $G$ we have that either $W_y$ is empty or equals $t^G\cap C_G(y)x$ for any $x\in W_y$. In particular, we have that $|W_y|\le |C_G(y)|$ for every $y\in G\setminus Z$. Therefore, this yields:
\begin{align*}
|W| = \sum_{y\in G} |W_y| & \le (1 + i(Z))|t^G| + \sum_{y\in G\setminus Z} |C_G(y)| \\ 
& = (1 + i(Z))|t^G| + \sum_{i=1}^r |y_i^G||C_G(y_i)|, 
\end{align*}
where $y_1,\ldots, y_r$ are the representatives of the non-central conjugacy classes of $G$. Thus $r=k(G)-|Z|$ and so
\begin{eqnarray}\label{eq1}
|W| \le (1 + i(Z)) \frac{|G|}{|C_G(t)|}  + (k(G) - |Z|) |G|.
\end{eqnarray}
On the other hand, observe that every element $x\in t^G$ inverts all elements of $[x,G]$, since $x^{-2}$ is central and so
\[
[x,g]^{x} = x^{-2} g^{-1} x g x = g^{-1} x^{-1} gx = [x,g]^{-1}
\] for any $g$ in $G$. As $|[x,G]|=|x^{-1}x^G|=|t^G|$ for $x\in t^G$, we then have that
\begin{eqnarray}\label{eq2}
\left|W \right| \ge \sum_{x\in t^G} |[x,G]| = |t^G|^2 = \frac{|G|^2}{|C_G(t)|^2}.
\end{eqnarray}
Hence, comparing (\ref{eq1}) and (\ref{eq2}) we get the desired equation.
\end{claimproof}

Now, let $x$ be an element of $G\setminus Z$ such that the order $|C_G(x)|$ is the maximum of all orders of centralizers for non-central elements, {\it i.e.}
\[
|C_G(x)| = \max \left\{ |C_G(y)| : y\in G\setminus Z \right\}.
\]Then, the class equation yields that
\[
|G| \ge |Z| + ( k(G) - |Z| )\frac{|G|}{|C_G(x)|}
\]
and so  $k(G) - |Z| < |C_G(x)|$, since certainly $|G|<|Z|+|G|$. Thus, we get that $k(G)-|Z|\le |C_G(x)|-1$. Combining this with the equation given by \nameref{Claim}, it follows that
\begin{align*}
|G| & \le (1 + i(Z)) |C_G(t)| +  ( |C_G(x)| - 1 ) |C_G(t)|^2  \\  
 & \le  |C_G(t)|^2 \left( |C_G(x)| - 1 + \frac{|Z|}{|C_G(t)| } \right),
\end{align*}
since $1+i(Z) \le |Z|$. This yields the desired inequality. 

\bibliographystyle{plain}

\end{document}